\documentstyle[12pt,amstex]{amsart}
\title{Algebraic geometry of Hopf-Galois extensions.}
\author{Dmitriy Rumynin}
\address{Department of Mathematics \& Statistics, 
LGRT, UMass, Amherst, MA, 01003.}
\email{rumynin@@math.umass.edu}
\subjclass{Primary 16W30; Secondary 14F05}
\date{September 24, 1997.}

\newtheorem{theorem}{Theorem}
\newtheorem{prop}[theorem]{Proposition}
\newtheorem{lemma}[theorem]{Lemma}

\newtheorem{cor}[theorem]{Corollary}

\newcommand{\Hom}{\mbox{Hom}}

\newcommand{\oo}{{\overline{O}}}
\newcommand{\ca}{{\mathfrak GAL}_{H, \mathcal S}}
\newcommand{\sss}{{\mathcal S}}
\newcommand{\uu}{{\mathcal U}}
\newcommand{\mo}{{\mathcal O}}
\newcommand{\ii}{\mbox{Id}}
\begin{document}
\maketitle

\begin{abstract}
We continue the investigation of Hopf-Galois extensions
with central invariants started in \cite{rum}. 
Our objective is not to imitate algebraic geometry using
Hopf-Galois extension but to understand their geometric 
properties.
\end{abstract}

Let $H$ be
a finite-dimensional Hopf algebra over a ground field
$\bf k$. Our main object of study  is an $H$-Galois extension
$U \supseteq O$
such that $O$ is a central subalgebra of $U$. Let us briefly
discuss geometric properties of the object. By Kreimer-Takeuchi
theorem the module $U_O$ is projective. Thus, it defines
a vector bundle of algebras on the spectrum of $O$
 by Serre theorem.
The fibers carry a structure of Frobenius algebra. A similar
structure was of interest to geometers for a while because
commutative Frobenius algebras naturally arise in the study
of symmetric Poisson brackets of hydro-dynamical type
 \cite{bal}. More recently, a concept of Frobenius
manifold was introduced \cite{hit}; it is a manifold
such that tangent spaces carry a structure of a commutative
Frobenius algebra which multiplication has a generating
function. 

Our set-up is different: we have a vector bundle
rather than the tangent bundle and 
our algebras are not necessarily
commutative. However, we have more structure involved:
a Hopf-Galois extension may be regarded as a ``quantum''
principal bundle \cite{sch}. 
(We should point out that the notion of a quantum
principal bundle in non-commutative geometry is
 more involved but, nevertheless,
quantum principal bundles with universal differential 
calculus are the same as Hopf-Galois
extensions \cite{dur, haj}.) 
If $H$ is commutative (i.e. an algebra of 
functions on a finite group
scheme $G$) then a commutative $H$-Galois
 extension $U \supseteq O$ 
is a $G$-principal bundle on the spectrum of $O$. 

Finally, we emphasize that centrality of invariants is
a crucial property for a geometrical treatment of Hopf-Galois
extensions. We will illustrate this claim throughout the paper.

We introduce Hopf-Galois extensions with central invariants
and discuss
examples in Section 1. The notion of inverse image of
a Hopf-Galois extension is discussed in Section 2.
A geometric object should become trivial if one looks at
it locally. We prove a weak version of this principle in Section 3.
A Hopf-Galois extension starts enjoying being
a geometrical object on an affine scheme 
in the next section where we show that it may be pasted
from local datum. The next three sections (5,6,7) are devoted
to investigation of various features a Hopf-Galois extension
carries. In section 8 we give a definition of Hopf-Galois
extension (or $H$-torsor) over any scheme.
and discuss the stack of $H$-torsors.

The author is grateful to  E. Markman, A. Mavlyutov,
  I. Mirkovic and  I. Taimanov
for the fruitful discussions.

\section{Introduction}

Let $H$ be a finite dimensional Hopf algebra over a ground field
{\bf k} from now on. One could extend a number of results
to the case of a ground ring  but it is beyond
our interest in the present paper. An associative algebra
with unity $U$ is an {\em $H$-comodule algebra}
if there is a map
$\rho :U \longrightarrow U \otimes H$, 
written as $\rho (x)=x_0 \otimes x_1$
such that $x_0 \varepsilon (x_1) =x$, 
$x_0 \otimes \Delta (x_1) = \rho (x_0) 
\otimes x_1$, and $\rho (xy)= \rho (x) \rho (y)$
 for all $x,y \in U$.
We use the oxymoronic variant of the Sweedler's $\Sigma$-notation
with $\Sigma$ eliminated. 
The subalgebra $O= U^H = \{ x \in U \mid  \rho (x)=x 
\otimes 1 \}$ is called {\em the subalgebra
of invariants}; one says that $U \supseteq O$ 
is {\em an $H$-extension}.
Throughout the paper we make various assumptions on $O$
as finitely generated, reduced (semiprime), or affine
(finitely generated semiprime).
An $H$-extension is {\em an extension with central invariants}
if $O$ is contained in the center of $U$. 
An $H$-extension is {\em Hopf-Galois}
(or specifically {\em $H$-Galois}) 
if the canonical map can:$U \otimes_O U 
\longrightarrow U \otimes H$
defined by can$(x \otimes_O y)=(x \otimes 1)(y_0 \otimes y_1)$
 is a bijection.
\cite{krt,mon}.

Let us discuss some examples of Hopf-Galois extensions.
A Galois extension of fields $K \supseteq {\bf k}$ with the Galois
group $G$
is Hopf-Galois for the Hopf algebra ${\bf k}G^{\ast}$. 
Geometrically, it means that the point over $K$ is a $G$-principal
bundle on the point over {\bf k}. 
The vector bundle which structure it carries
has the {\bf k}-vector space $K$ as a fiber at the point.

However, a Hopf-Galois extension of fields 
is not necessarily 
Galois. A purely inseparable extension may be Hopf-Galois
for the dual of a restricted universal enveloping
algebra \cite{mon}.
Moreover, a finite separable not normal field extension
may be Hopf-Galois as well \cite{grp}. 
From the geometric prospective 
it means that the point over $K$ 
is a principal
bundle on the point over {\bf k} for some 
finite non-reduced group scheme.

A criterion  for a finite separable
field extension to be Hopf-Galois was proven in \cite{grp}.
One should consider the normalization 
$\tilde{K}$ of $K$ over {\bf k}.
Let $S(M)$ be the symmetric group of the set
 $M= \mbox{Gal}( \tilde{K} \mid {\bf k})/
\mbox{Gal}( \tilde{K} \mid K)$. There is a natural 
embedding $ \mbox{Gal}( \tilde{K} 
\mid {\bf k})
\hookrightarrow S(M)$. The extension is Hopf-Galois
if and only if there exists a subgroup $N \subseteq S(M)$, 
normalized by 
$ \mbox{Gal}( \tilde{K} \mid {\bf k})$,
such that the left $N$-set $M$ is isomorphic 
to the left regular $N$-set $N$.

Interesting Hopf-Galois extensions may be 
constructed as crossed products.
More precisely, crossed products constitute all 
cleft Hopf-Galois extensions.
An inquisitive reader should consult \cite{mon} 
for a nice account
of the subject. We do not need the general 
construction. We restrict
our attention to twisted products, crossed
 products with the trivial
action, because only they produce central invariants.
Let $R$ be a commutative ${\bf k}$-algebra and 
$\sigma : H \otimes H \longrightarrow R$
be a linear map. Let $R_{\sigma}[H]$ be a vector space
$R \otimes H$ with an algebra structure defined by the formula
\begin{equation} a \otimes g \cdot b \otimes h =  
ab \sigma (h_1,g_1) \otimes h_2g_2
\ , \ \ \ a,b \in R , g,h \in H \label{mult} \end{equation}
It may not be a structure of an associative algebra in general.
The following lemma is straightforward \cite[Lemma 7.12]{mon}.

\begin{lemma}
The formula~\ref{mult} defines a structure
 of an associative algebra
with an identity element $1 \otimes 1$ if and only if
for each $g,h,t \in H$
$$ \sigma (h_1 \otimes g_1) \sigma (h_2 g_2 \otimes t)=  
\sigma (g_1 \otimes t_1) \sigma (h \otimes g_2t_2)$$
$$ \sigma (h \otimes 1) = \sigma (1 \otimes h)
 = \varepsilon (h) 1$$
\end{lemma}

If $\sigma$ is invertible with respect 
to the convolution \cite{mon}
and satisfies the conditions of Lemma 1 
we call it {\em a cocycle},
although the corresponding cochain complex has been constructed
only for cocommutative $H$ \cite{swe}.
We use the term {\em twisted product} for $R_{\sigma}[H]$ 
if $\sigma$ is a cocycle.
The natural question is what is happening 
if the map $\sigma$ is not
invertible but nevertheless  formula~\ref{mult} defines 
a structure of an associative algebra.
It may happen if there is a family of 
cocycles $\sigma_t$ converging
in an appropriate sense to a linear map 
$\sigma_0$; the latter inherits
the identities of Lemma 1 but can fail 
to be convolution  invertible.
In other words, it is a closed condition 
to define an associative algebra but 
it is an open one to be invertible. 
For such a ``non-invertible cocycle'' one 
still obtains an $H$-comodule
algebra which is not Hopf-Galois. However, 
an important piece of 
the structure theory
breaks down in this situation.

\begin{lemma} 
Assume that $U \supseteq O$ is an $H$-extension.
The following statements are equivalent

1. $U$ is isomorphic to a twisted product 
$O_{\sigma} [H]$ as an 
$H$-comodule algebra.

2. There exists a convolution invertible 
right $H$-comodule map
$\gamma : H \longrightarrow U$ and 
the subalgebra $O$ is central.

3. There is a linear map from $U$ to 
$O \otimes H$ conducting 
an isomorphism
of left $O$-modules and  right $H$-comodules 
and $O$ is central.

4. $U$ is a free $O \otimes H^{\ast}$-module of rank 1
and $U^H$ is central.
\end{lemma}
{\bf Proof.} The equivalence of 1 through 3 
follows from
Theorems 7.2.2 and 8.2.4 of \cite{mon} and 
Lemma 11 of \cite{rum}.
4 follows from 1 because $1 \# \Lambda$ for 
a non-zero left
integral $\Lambda$ is clearly a basis of 
the left $O \otimes H^{\ast}$-module
$O_{\sigma}[H]$. Finally, 4 implies 3 since 
$O \otimes H$ is a free left
$O \otimes H^{\ast}$-module of rank 1. Thus, 
there exists an 
$O \otimes H^{\ast}$-module isomorphism between 
$U$ and $O \otimes H$
which must be an isomorphism of left $O$-modules 
and right $H$-comodules.
$\Box$

Now it is clear why only crossed products with the trivial
action are of interest to us: the action must 
be trivial to produce central invariants.

Let us look at an example.
We define the following cocycle on ${\Bbb C}{\Bbb Z}_n$.
Denoting a generator of ${\Bbb Z}_n$ by $t$ we set for each
$0 \leq k,m < n$
\begin{equation} \sigma (t^k \otimes t^m)= 
\left\{ \begin{array}{ll}
 1 & \mbox{ if } k+m < n \\ 
q & \mbox{ if } k+m \geq n  \end{array} \right.  
\label{mu} \end{equation}
It is a cocycle if and only if $q \in O$ 
is an invertible
element. The first choice may be  $q 
\in O={\Bbb C}[q,q^{-1}]$. 
The twisted product 
${\Bbb C}[q,q^{-1}]_{\sigma}[{\Bbb C}{\Bbb Z}_n]$
is the small quantum cohomology ring 
of the projective space ${\Bbb CP}^{n-1}$.
The multiplication coefficients of it in an appropriate basis
are 3-point Gromov-Witten 
invariants. One should consult \cite{mcds} for a good account
of the subject. Another choice one can make is
$q=0 \in {\Bbb C}$. The corresponding ring
${\Bbb C}_{\sigma}[{\Bbb C}{\Bbb Z}_n]$ is the de Rham
cohomology ring of ${\Bbb CP}^{n-1}$.
It is an example of a ``non-invertible twisted product'':
it may be thought of as a limit of twisted products depending
on $q \in {\Bbb C} \setminus \{ 0\}$ as $q \rightarrow 0$. 
One can make another choice of a ring: let $q=e^{x_1} \in
C^{\infty}({ \Bbb C}^n)$. After identification
$t^{k-1}= \frac{\partial}{\partial x_k}$ one obtains a structure
of Frobenius manifold on ${\Bbb C}^n$ \cite{hit}. 

This example 
 demonstrates that one can encounter an interesting geometry
studying Hopf-Galois extensions with central invariants. From 
now on $U \supseteq O$ is an $H$-Galois
extension with central invariants unless a contrary
is specified. 
Let $\Lambda$ be a non-zero left integral of $H^{\ast}$. 
By the definition
$\Lambda \cdot x = x_0 \Lambda (x_1)$.
The following facts were first proven in 
\cite[1.7, 1.9, 1.10]{krt}. One may find an interesting
discussion concerning the last two properties in \cite{co1}.
\begin{lemma} The following statements hold under our assumptions.

1. $U_O$ is a projective finitely generated module.

2. $U_O$ is faithfully flat.

3. There exists an $O$ submodule $P \subseteq U$ such that
$U=O \oplus P$.

4. There exists an element $x \in U$ such that 
$\Lambda \cdot x =1$.
\end{lemma}

We can now explain why centrality of invariants is essential.
One may hope, for instance, that $O$ being commutative suffices
to build an interesting geometry. Indeed, there is a one-to-one
correspondence
between algebraic vector bundles on the spectrum of $O$
and finitely generated projective modules by Serre theorem 
\cite[Corollary 2.4.50]{ser}.
If $O$ is not central then left and right actions are different.
Thus, we get  {\em
a bi-bundle} (i.e. a space carrying two commuting 
structures of a vector
bundle) rather than a vector bundle. Geometry 
of bi-bundles is less exploded 
than that of vector bundles. Thus, if we want 
to remain in the realm
of vector bundles we should require invariants 
being central. In further sections
we will see more properties which do not work for Hopf-Galois
extensions with non-central invariants.

We finish the introduction with the lemma which
is  handy if one wants to test whether an $H$-extension
is Hopf-Galois. It is a slight 
generalization of \cite[Lemma 1.3]{krt}.

\begin{lemma} If there exists a subalgebra $A \subseteq O$ such 
that the map
$\Gamma : U \otimes_A U \longrightarrow U \otimes H$ defined by
$\Gamma (x \otimes_O y)=(x \otimes 1)(y_0 \otimes y_1)$ 
is onto then $U \supseteq O$ is $H$-Galois.
\end{lemma}
{\bf Proof.} The map $\Gamma$ factors 
through the canonical map.

\begin{picture}(150,100)(0,0)
\put(10,80){$U \otimes_O U$}
\put(90,80){$U \otimes H$}
\put(53,83){\vector(1,0){30}}
\put(59,86){can}
\put(10,20){$U \otimes_A U$}
\put(25,35){\vector(0,1){35}}
\put(40,35){\vector(3,2){52}}
\put(61,40){$\Gamma$}
\end{picture}

Thus, the canonical map is onto. 
It must be also one-to-one by \cite[Theorem 8.3.1.]{mon}.

$\Box$

We should notice that $\Gamma$ being 
one-to-one in the last lemma does not
imply $A=O$. A proper embedding $A \hookrightarrow O$ 
may be an epimorphism
in the category of rings making the tensor products 
$U \otimes_O U$
and $U \otimes_A U$ indistinguishable.

\section{Pull-back.}

This section is devoted to the crucial geometric property
of Hopf-Galois extensions with central invariants, the pull-back. 
Intuitively, given a map of schemes $X \longrightarrow Y$ and
a geometrical object on $Y$, we should be able to induce a similar
object using some kind of the fiber product. So far, we can work
only on the level of rings. 
The following lemma provides the existence
of pull-back on the level of rings. 
This property fails for Hopf-Galois
extensions with non-central invariants.
\begin{lemma} Let $\phi :O \rightarrow R$ be a map of commutative
algebras. Then $\tilde{U}=U \otimes_O R \supseteq R$ 
is an $H$-Galois
extension with central invariants.
\end{lemma}
{\bf Proof.}
The natural map $R \rightarrow \tilde{U}$
is one-to-one by Lemma 3. Indeed, $U = O \oplus P$ 
as $O$-modules
and a split extension is  pure. In down-to-earth terms,
$\tilde{U}=(O \oplus P) \otimes_O R= 
O \otimes_O R \oplus P \otimes_O R
\supseteq O \otimes_O R \cong R$.

The coaction of $H$ on $\tilde{U}$ is given by
$(u \otimes_O r)_0 \otimes (u \otimes_O r)_1= 
(u \otimes_O r_0) \otimes r_1$.
Let $S$ be the subalgebra of invariants $\tilde{U}^H$. 
It is clear that
$S \supseteq R$. Let us consider the diagram.

\begin{picture}(150,100)(0,0)
\put(10,80){$\tilde{U} \otimes_S \tilde{U}$}
\put(90,80){$\tilde{U} \otimes H \cong S \otimes_O (U \otimes H)$}
\put(50,83){\vector(1,0){30}}
\put(55,87){can}
\put(10,20){$\tilde{U} \otimes_O \tilde{U}$}
\put(25,35){\vector(0,1){35}}
\put(40,35){\vector(3,2){54}}
\put(65,43){$\Gamma$}
\end{picture}

The map $\Gamma$ is onto since 
the original extension is Hopf-Galois.
By Lemma 4 the extension $\tilde{U} \supseteq S$ is $H$-Galois.

It remains to show that $S=R$.
Let $x \in U$ be
an element such that $\Lambda \cdot x = 1$ 
which existence is provided
by Lemma 3. It is clear that
$\Lambda \cdot (x \otimes_O 1) = 
1 \otimes_O 1 = 1 \in R$. It implies
$\Lambda \cdot  \tilde{U} = \tilde{U}^H$. 
But $\Lambda \cdot  \tilde{U} 
\subseteq \Lambda \cdot  ((U \otimes_O 1)R) \subseteq
R ( \Lambda \cdot  U \otimes_O 1) \subseteq 
R ( ( \Lambda \cdot  U) \otimes_O 1) \subseteq R(O \otimes_O 1) 
\subseteq R$ 
 $\ \ \ \Box$

The module $U_O$ defines a vector bundle 
of algebras $\tilde{U}$ 
on Spec$\: O$ \cite{rum}. Let us denote by 
$O_{(\chi)}$ the local ring of
algebraic functions at 
$\chi \in \mbox{Spec} \: O$, $m_{\chi}$ its
maximal ideal, and 
$K_\chi$ the quotient field $O_{(\chi)} / m_\chi$.
We remind the main characters of the play.

{\em The germ algebra} may be defined by
\begin{equation} U_{( \chi )} =
U \otimes_{O} O_{(\chi )} \supseteq O_{(\chi)}.
\label{2} \end{equation}
Here is another way to describe it.
$S= O \setminus \chi$ is a multiplicative subset of $U$. 
The generalized algebra of quotients $US^{-1}$ 
is isomorphic to $ U_{(\chi)}$.
{\em The fiber algebra} is
\begin{equation} U_\chi =
U \otimes_{O} K_{\chi } \supseteq K_{\chi}.
\label{3} \end{equation}
One may introduce {\em the $n$-jet algebra} 
generalizing the germs
and the fibers. It is important 
for the study of indecomposable
modules in the spirit of \cite{rud}.
\begin{equation} U_\chi^{(n)} =
U \otimes_{O} O_{(\chi )}/m_{\chi}^n 
\subseteq O_{(\chi)}/m_{\chi}^n.
\label{4} \end{equation}
It is obvious that
$U_{\chi}^{(n)} \cong U_{(\chi )}/ U_{(\chi )}m_{\chi}^n$.
 Conventionally, $m_{\chi}^\infty =0$. 
Thus, the algebra of $\infty$-jets 
$U^{( \infty )}_{\chi}$
is the same as the germ algebra $U_{(\chi )}$.
Another interesting choice is $n=1$: 
$U_{\chi}^{(1)}$ is the fiber 
algebra $U_{\chi}$ at the point $\chi$. 
The last algebra we want to introduce is the restriction
to the closure of a point $\chi$:
$$U_{[\chi ]}=U / U \chi \supseteq O/ \chi .$$
If the point $\chi$ is closed then 
$U_{[\chi]}$ is isomorphic to $U_ \chi $.

Now we can prove the main theorem of the section
improving Theorems 15 and 16 of \cite{rum}.

\begin{theorem} 
$U_{[\chi]}^{H} \supseteq O/ \chi$ and
$U_{\chi}^{(n)} \supseteq O_{\chi}/m_{\chi}^n$
are $H$-Galois extensions for each $\chi$ and $n$
(including $n= \infty$).
If the extension $U \supseteq O$ is cleft then 
so are the jet extensions and the restrictions.
\end{theorem}

{\bf Proof.}
The first statement is an immediate corollary of Lemma 5. 
To prove the second statement we assume that
$\gamma :H \longrightarrow U$ is a splitting map.
Then the  composition  
$H \stackrel{\gamma}{\longrightarrow} U \longrightarrow
U_{\chi}^{(n)}$
is a splitting map for the jet extension  
$U_{\chi}^{(n)} \supseteq O_{( \chi )}
/m_{\chi}^n$.
Similarly, the  composition  
$H \stackrel{\gamma}{\longrightarrow} U \longrightarrow
U/U \chi$
is a splitting map for the restriction extension  
$U/U \chi \supseteq O/ \chi$.
$\Box$

This theorem is not the best statement one can get. 
We will prove
in the next section that the jet algebras are always cleft for
$n < \infty$ and generically cleft for 
$n=\infty$ under mild restrictions
on $O$. Let us give one
more definition. Let ${\frak GAL}_H(O)$ be a set of isomorphism
classes of $H$-Galois extensions $U \supseteq O$. We understand
an isomorphism of $H$-comodule algebras by an isomorphism.
Following \cite{rum}, we denote 
the subset of isomorphism
classes of cleft $H$-Galois extensions by ${\frak Gal}_H$.
The following theorem is an immediate consequence of Lemma 5.
It improves \cite[Theorem 14]{rum}.

\begin{theorem}
${\frak GAL}_H$ is a covariant functor from the category of 
commutative {\bf k}-algebras 
to the category of sets and ${\frak Gal}_H$ is its subfunctor.
\end{theorem}

\section{Local structure.}

Usually if one looks at a geometrical 
object locally it seems trivial. 
A trivial $H$-Galois extension is 
the tensor product $O \otimes H \supseteq O$.
It is false that an $H$-Galois extension 
is locally trivial: if $U \supseteq O$
is a proper Galois extensions of fields 
with Galois group $G$ then $O$ is 
already localized, nevertheless $U$ 
is not isomorphic to $O \otimes {\bf k}G^\ast$.
However, one could expect that a Hopf-Galois 
extension is locally cleft.
We don't know whether it holds. 
The following question is interesting from this
prospective. 
In this section we prove  a weaker version that
a Hopf-Galois extension is generically locally cleft.

{\bf Question.} Let $U \supseteq O$ be an 
$H$-Galois with central $O$ 
and finite-dimensional
$H$. Assume that $O$ is a local algebra. 
Is the extension cleft?

First, we want to understand what is happening 
if the vector bundle
defined by a Hopf-Galois extension is trivial. If we managed
to show that this implied that the Hopf-Galois extension were cleft
it would be very nice. Indeed, any vector bundle is locally trivial.
The following fact is the best we can do. It first appeared in
\cite[Lemma 2.9]{bla}. We should mention that trivial bundles
correspond to free modules on the algebraic side.

\begin{lemma} 
If $U_O$ is a free module of rank $n$ then 
$U^n$ is a free left $O \otimes H^{\ast}$-module
of rank $n$.
\end{lemma}
{\bf Proof.}
Using the canonical map we get a chain of isomorphisms:
$$ U^n \cong U \otimes_O U \cong U \otimes 
H \cong U \otimes_O (O \otimes H)
\cong (O \otimes H)^n  \cong (O \otimes H^{\ast})^n 
\ \ \ \ \Box$$

Thus, $U$ is a projective module of rank 1 
if the rank is well-defined for
projective $O \otimes H^{\ast}$-modules. 
The next idea is to put
some constraints on $O$ and $H$ to ensure that any projective 
$O \otimes H^{\ast}$-module of rank 1 is free. A
 reader may supplement
the next corollary with other similar statements 
providing cleftness.

\begin{cor}  $U$ is cleft if one of 
the following conditions holds:

1. $O$ is artiniian and $U_O$ is free;

2. $O$ is a local artiniian ring;

3. $O$ is a field;

4. $O \otimes H^{\ast}$ is self-injective and $U_O$ is free;

5. $\mbox{\em Spec} \, O$ is a toric variety 
and $O \otimes H^{\ast}$
is self-injective. 
\end{cor}
{\bf Proof.} Each of the conditions 1-5 implies 
the module $U_O$ is free.
It is explicitly stated in 1 and 4. Any finitely
generated projective module over local ring or field is free.
Algebra of functions on a toric variety 
is a semigroup algebra of
a normal monoid.  
A finitely generated projective module 
is free over such an algebra
according to the Anderson conjecture 
proved in \cite[Theorem 2.1]{gub}.

Using Lemma 8, it suffices to  understand 
some kind of uniqueness
of decomposition theorem in each of the cases.
In 1, 2, and 3 one can use 
the Krull-Remark-Schmidt theorem \cite[5.18.11]{fai}.
Indeed, $O \otimes H^{\ast}$ is
artiniian if so is $O$: since $H$ is finite 
dimensional $O \otimes H^\ast$
is even an artiniian $O$-module.
Thus the $O \otimes H^{\ast}$-module
$U^n$ must have a finite length which 
implies the uniqueness
of the decomposition.

We can use a stronger version of the 
Krull-Remark-Schmidt theorem
in 4 and 5 \cite[5.19.18]{fai}. 
$U^n$ is injective $O \otimes H^\ast$-module
in this
case; so are its direct summands. The endomorphism ring
of an injective indecomposable module 
is local and we can use 
the Krull-Remark-Schmidt
theorem.
$\Box$ 
 
Now we can state and prove the local trivialization theorem.

\begin{theorem}
The germ extension is  cleft 
for generic $\chi \in \mbox{Spec} \; O$.
If $O$ is finitely generated  then the $n$-jet extension
is cleft for every $\chi \in \mbox{Spec} \; O$ and $n < \infty$. 
\end{theorem}
{\bf Proof.} 
The second statement follows from Corollary 9.2 since
the $n$-jets of functions $O_{(\chi )}/m_{\chi}^n$ is a local
artiniian ring (even finite dimensional) in the case of
finitely generated $O$.

The first statement can be reduced 
to a component of the spectrum:
let $I$ be a minimal prime ideal of 
$O$. We replace the extension
$U \supseteq O$ with $U \otimes_OO/I \supseteq O/I$.
Thus, we can assume $O$ is a domain 
without loss of generality.
This argument becomes sloppy unless $O$ is finitely
generated. Otherwise, no minimal prime ideal can exist.
However, the term generic is unclear 
in this case. Thus, we can make
up the definition of a generic point and 
think that the theorem still
holds. For example, generic could 
mean generic in closure of each point:
it entitles us to replace $O$ with $O/I$ 
where $I$ is an arbitrary
prime ideal.

We should use the field fractions 
$Q(O)$ to comprehend the first statement. 
Geometrically, it is a field
of rational functions on the spectrum. By Corollary 9.3  the
extension $U \otimes_O Q(O) \supseteq Q(O)$ is cleft.
Let $\gamma :H \longrightarrow U \otimes_O Q(O)$ be
a splitting map. The algebra $U \otimes_O Q(O)$ can be realized
as a generalized ring of quotients of $U$ by the multiplicative
set $O \setminus \{ 0 \}$. Let $h_i$ be a base of $H$.
Let $\mathcal V$ be the complement of zeroes 
of the denominators of $\gamma (h_i)$.
It is obvious that $U_{(\chi )} \supseteq O_{(\chi )}$ is cleft
for each $\chi \in {\mathcal V}$  $\ \ \ \Box$

\section{Pasting.}

We do algebraic geometry of $H$-Galois 
extensions on affine schemes
on the level of rings till Section 8 
where we sheafify the enterprise.
After Section 2, given an $H$-Galois 
extension of an affine scheme $X$
and an open cover $\cup V_i$ of $X$, 
we can restrict the extension
to each set $V_i$. 
Now we are interested in the inverse problem:
given $H$-Galois extensions on each $V_i$, 
compatible on intersections,
we want to paste an $H$-Galois extension on $X$.

We are working in the flat topology since it 
is sufficiently general for a number
of applications. {\em A covering} of an algebra 
$O$ is a finite set of $O$-algebras $O_i$ 
such that the product $\prod_i O_i$ 
(denoted $\overline{O}$ from now on) is a finitely
presented $O$-algebra which is a faithfully flat $O$-module. 
The last condition provides
that the natural map of spectra 
$\coprod_i \mbox{Spec} \, O_i \rightarrow \mbox{Spec} \, O$
is onto.

The standard example of a covering can be 
constructed using the partition of unity.
Given a finite set of elements $\{ f_i \}$ 
of $O$ generating the trivial ideal $O$
(i.e. $\sum_ix_if_i=1$ for some $x_i$), 
we set $O_i=Of_i^{-1}$.
It is standard to check that we get a covering this way. 
We also point out that
partitions of unity are the same 
as coverings in Zariski topology. 

We introduce another piece of notation: 
given a covering $O_i$ we define 
$O_{ij}=O_i \otimes_OO_j$
and $O_{ijk}=O_i \otimes_OO_j \otimes_OO_k$. 
These algebras should be thought of as double and triple
intersections in the covering.

We consider a covering $O_i$ of $O$ and 
a collection of $H$-Galois extensions $U_i \supseteq O_i$.
We assume there are  isomorphisms
of $H$-comodule algebras 
$\phi_{ij}:U_i \otimes_{O_i}O_{ij} 
\rightarrow U_j \otimes_{O_j}O_{ij}$.
The datum $(O_i,U_i, \phi_{ij})$ 
is called {\em an $H$-structure}
on the algebra $O$
if $\phi_{ij}$, restricted to $O_{ij}$,
is the identity map and  
the cocycle condition is satisfied:
$(\phi_{jk} \otimes_{O_{jk}}\ii_{O_{ijk}}) \circ 
(\phi_{ij}\otimes_{O_{ij}}\ii_{O_{ijk}})=
\phi_{ik}\otimes_{O_{ik}}\ii_{O_{ijk}}$
for each $i,j,k$.

\begin{theorem} For each $H$-structure on an algebra $O$
there exists a unique up to an isomorphism 
$H$-Galois extension with central invariants 
$U \supseteq O$ such that
$U \otimes_O O_i \cong U_i$.
\end{theorem}
This theorem allows to glue Hopf-Galois 
extensions from local data. 
However, this locality property is 
insufficient to conclude ${\frak GAL}_H$ is a sheaf
in the flat topology in the category of 
$\bf k$-algebras because there can be non-trivial
automorphisms. For instance, if $G$ is 
a group with non-trivial center $Z$ and $U \supseteq O$
is a $G$-Galois field extension then $Z$ 
acts by ${\bf k}G^\ast$-comodule algebra
automorphisms. Nevertheless, ${\frak GAL}_H$ 
can made into a stack which will be discussed in Section 8.

A natural question is to understand what kind 
of local property the subfunctor ${\frak Gal}_H$
inherits. Given an $H$-structure with split extensions $U_i$,
 the global extension $U \supseteq O$,
whose existence is ensured by Theorem 11, 
is not necessarily split.
Indeed, let $\gamma_i :H \rightarrow U_i$ 
be a bunch of splittings (i.e. convolution-invertible right
$H$-comodule maps); one has to glue a splitting 
$\gamma :H \rightarrow U$. It can be done 
if $\gamma_i$ agree on ``intersections'' $U_{ij}$
but it is unclear how to do it without any 
restrictions on behavior on intersections.

The rest of the section is devoted to the proof 
of Theorem 11. The major tool is
the following lemma \cite[1.1.3.14]{dem}. 
If $f:A \rightarrow B$ is a morphism of rings
then a linear map $h:M \rightarrow N$ from 
an $A$-module $M$ to a $B$-module $N$
is called {\em adapted to $f$} if 
$f(am)=h(a)f(m)$ for each $a \in A, m \in M$.

\begin{lemma} Let $R$ be a commutative ring and 
$B$ be a faithfully flat commutative $R$-algebra.
Let us consider modules $J$, $K$, $L$ over $B$, 
$B \otimes_RB$, and 
$B \otimes_RB \otimes_RB$ with homomorphisms
\begin{center}
\begin{picture}(100,60)(0,0)
\put(0,20){$J$}
\put(10,20){\vector(1,0){19}}
\put(12,9){$u_1$}
\put(10,28){\vector(1,0){19}}
\put(12,35){$u_0$}
\put(30,20){$K$}
\put(43,14){\vector(1,0){19}}
\put(46,2){$u_2^\prime$}
\put(43,19){\vector(1,0){19}}
\put(46,24){$u_1^\prime$}
\put(43,35){$\vector(1,0){19}$}
\put(46,41){$u_0^\prime$}
\put(70,20){$L$}
\end{picture}
\end{center}
adapted to corresponding ring homomorphisms
\begin{center}
\begin{picture}(140,60)(0,0)
\put(0,20){$B$}
\put(11,20){\vector(1,0){18}}
\put(12,9){$d_1$}
\put(11,26){\vector(1,0){18}}
\put(12,31){$d_0$}
\put(32,20){$B \otimes_RB$}
\put(75,14){\vector(1,0){19}}
\put(78,3){$d_2^\prime$}
\put(75,19){\vector(1,0){19}}
\put(78,25){$d_1^\prime$}
\put(75,37){\vector(1,0){19}}
\put(78,43){$d_0^\prime$}
\put(99,20){$B \otimes_RB \otimes_RB$}
\end{picture}
\end{center}
such that $d_0(b)=b \otimes 1, d_1(b) = 1 \otimes b, 
d_0^\prime (b \otimes c)= b \otimes c \otimes 1,
d_1^\prime (b \otimes c)= b \otimes 1 \otimes c$, 
and $d_2^\prime (b \otimes c)= 1 \otimes b \otimes c$
for each $b,c \in B$. We also assume that 
$u_0^\prime u_0=u_1^\prime u_0, 
u_0^\prime u_1=u_2^\prime u_0$,
and $ u_1^\prime u_1=u_2^\prime u_1$. If 
$I= \mbox{Ker}(u_0,u_1)$ then the embedding 
$I \hookrightarrow J$
induces an isomorphism $I \otimes_RB \rightarrow J$.
\end{lemma}

{\bf Proof.}
Let us first look at the special case $O= \prod O_i$. 
It is easy to see that $U= \prod_i U_i$
satisfies the conditions of the theorem. 
If $\check{U}$ is another extension satisfying
the conditions of the theorem then 
$\check{U} \cong \check{U} \otimes_OO \cong
\prod_i \check{U}\otimes_OO_i \cong \prod_iU_i \cong U$. 
It is easy to see that all isomorphisms
are those of $H$-comodule algebras. 
We have just established an important property of
${\frak GAL}_H$ that for finite products
$${\frak GAL}_H (\prod O_i)= \prod_i {\frak GAL}_H (O_i).$$

Now we want to use Lemma 12 with $R=O$, 
$B= \overline{O} = \prod_i O_i$, $J=\prod_iU_i$,
$K=J \otimes_{\oo}^2 (\oo \otimes_O \oo )$ where 
2 means that the tensor product is taken
using the embedding $d_1$ to the second factor 
$b \mapsto 1 \otimes b$, and, finally,
$L=J \otimes_{\oo}^3 (\oo \otimes_O \oo \otimes_O \oo )$ 
given by the embedding to the third factor 
$b \mapsto 1 \otimes 1 \otimes b$. We also need some maps, 
three of which are easy to describe:
$$u_1(x)=x \otimes 1\otimes 1, \ u_1^\prime 
(x \otimes a \otimes b)
=x \otimes a \otimes 1 \otimes b, \ 
u_2^\prime (x \otimes a \otimes b)=
x \otimes 1 \otimes a \otimes b.$$
It is straightforward to see that $u_1$, $u_1^\prime$, 
and $u_2^\prime$ are adapted to
$d_1$, $d_1^\prime$, and $d_2^\prime$ correspondently.
Let us define
$$\phi = \prod_{i,j}\phi_{ij}:J {\otimes^1}_{\oo} 
(\oo \otimes_O\oo )\rightarrow
J {\otimes^2}_{\oo} (\oo \otimes_O\oo ).$$
$J \supseteq \oo$ is an $H$-Galois extension 
because the theorem
has already been proved for direct products. 
$\phi$ is an isomorphism
of $H$-Galois extensions. We can define
$$u_0(x)=\phi (x \otimes^1 (1 \otimes 1)).$$
Finally we need an isomorphism
$$\prod_{i,j,k}\phi_{ik}\otimes_{O_{ik}} \ii_{O_{ijk}}:
J\otimes^1_{\oo} (\oo \otimes_O \oo \otimes_O 
\oo) \rightarrow
J\otimes^3_{\oo} (\oo \otimes_O \oo \otimes_O \oo)$$
to define 
$$u_0^\prime (x \otimes a \otimes b)=
\prod_{i,j,k}\phi_{ik}\otimes_{O_{ik}} \ii_{O_{ijk}}
(\phi^{-1}(x \otimes a \otimes b) \otimes 1).$$
It is straightforward to check that 
$u_0$ and $u_0^\prime$ are adapted to $d_0$
and $d_0^\prime$ correspondently.
We have to compute the condition
for compositions to apply Lemma 12.
It is easy to see that
$$u_1^\prime (u_1(x))= x \otimes (1 \otimes 1 \otimes 1)=
u_2^\prime (u_1(x)).$$
Assuming $\phi (x \otimes 1 \otimes 1)=
\sum_t y_t \otimes a_t \otimes b_t$, we  compute that
$$u_0^\prime (u_0(x))= \sum_t y_t 
\otimes (b_t \otimes 1 \otimes c_t)=
u_1^\prime (u_0(x)).$$
Finally,
$$u_2^\prime (u_0(x))= \sum_t
y_t \otimes (1 \otimes a_t \otimes b_t).$$
Assuming  $\phi^{-1} (x \otimes 1 \otimes 1)=
\sum_t z_t \otimes c_t \otimes d_t$, we can write
$$u_0^\prime (u_1(x))=
\prod_{i,j,k}\phi_{ik}\otimes_{O_{ik}} \ii_{O_{ijk}}
(\sum_tz_t \otimes (c_t \otimes d_t \otimes 1)).$$
The results of the last 
two calculations because of the cocycle condition.

Now we define $U= \mbox{Ker}(u_0,u_1)=\{x \in J \mid \phi 
(x \otimes (1 \otimes 1))=x \otimes (1 \otimes 1) \}$.
$U$ is a subalgebra of $J$, closed under $H^\ast$-action,
since $\phi$ is a map of $H$-comodule algebras.
$U^H=U \cap J^H= U \cap \oo=\{x \in \oo \mid \phi 
(x \otimes (1 \otimes 1))=x \otimes (1 \otimes 1) \}=O$.
Finally, $\mbox{can}_J=\mbox{can}_U \otimes_O \ii_{\oo}$
modulo natural identifications 
$(U \otimes H) \otimes_O \oo \cong (U \otimes_O \oo )
\otimes H \cong J \otimes H$ and
$(U \otimes_O U) \otimes_O \oo \cong
(U \otimes_O \oo) \otimes_\oo (U \otimes_O \oo )
\cong J \otimes_\oo J$.
The second identification works through the map
$(a \otimes x) \otimes (b \otimes y) \mapsto
(a \otimes b) \otimes xy$.
Since $\oo_O$ is faithfully flat and $\mbox{can}_J$
is an isomorphism the map $\mbox{can}_U$ is
an isomorphism too and $U \supseteq O$ is $H$-Galois. 

Let $\check{U}$ be another extension with the given properties.
Then $\check{U} \otimes_O \oo \cong \prod_i 
\check{U} \otimes_OO_i \cong \prod_j U_i \cong J$.
Thus, there is the natural map $\check{U} \rightarrow J$
which factors through $U$ because of the universal
property of the kernel of a pair. 
The map $\check{U} \rightarrow U$
is of projective $O$-modules
which is a local isomorphism (on the covering $O_i$).
Thus, its kernel and cokernel have trivial support 
and the map in an isomorphism. 
Noticing that it is a map of $H$-comodule algebras
finishes the proof. 
$\ \ \Box$
 
The argument about the canonical maps fails 
for Hopf-Galois extensions
which invariants are not central because of the lack
of the identifications.

\section{Frobenius form.}

A Frobenius manifold has a Riemannian metric as a part of its
structure. Similarly, Hopf-Galois extensions carry a canonical
(up to a scalar) non-degenerate associative 
form which can  fail
to be symmetric.
Chosen $\Lambda$, a non-zero left integral
of $H^{\ast}$,
we construct an $O$-bilinear
form $\langle , \rangle :U \times U \longrightarrow O$ by
$\langle x,y \rangle = \Lambda \cdot (xy) = 
x_0y_0 \Lambda (x_1y_1)$ for each $x,y \in U$.

\begin{lemma} \cite[1.7.5]{krt}
The for $\langle , \rangle$ is non-degenerate.
\end{lemma}

Non-degeneracy means that the map $\zeta : U \longrightarrow
\Hom (U_O, O)$ given by 
$\zeta (u)(v)= \langle u,v \rangle$ is an isomorphism. 
The following corollary is immediate from
Lemma 13, Theorem 6, and the definition of a Frobenius algebra.

\begin{cor}
The algebras $U_{\chi}$ are Frobenius $K_{\chi}$-algebras.
\end{cor}

The question when this form is symmetric is quite interesting.
The following theorem is a sufficient condition for 
being symmetrical.
It is a slight generalization of \cite[Theorem 17]{rum}.
\begin{theorem}
Let $H$ be unimodular with the antipode 
of order 2. The form $\langle , \rangle$ is symmetric
if one of the following holds:

1. $U \supseteq O$ is cleft.

2. $O$ is semiprime.
\end{theorem}
{\bf Proof.}
Let us treat the cleft case first. We assume that $U$ is isomorphic
to a twisted product $O_{\sigma}[H]$.
We notice that
$h_1 \Lambda (h_2)= \Lambda (h) 1$
for each $h \in H$.
Indeed, the equality 
$ \alpha ( h_1 \Lambda (h_2))= \alpha ( \Lambda (h) 1)$
holds for each $\alpha \in H^{\ast}$ since
$ \alpha ( h_1 \Lambda (h_2))=  \alpha ( h_1) \Lambda (h_2)= 
\alpha  \Lambda (h) = \alpha ( 1) \Lambda (h)$.
It was proved in \cite{obe} that $H$ is unimodular
with the antipode of order 2 if and only if
$\Lambda (xy) = \Lambda (yx)$ for each $x,y \in H$.
Therefore, $\langle a \otimes h , 
b \otimes g \rangle = abh_1g_1 \Lambda (h_2g_2)=
ab \Lambda(hg) = ba \Lambda (gh)= \langle b \otimes g , 
a \otimes h \rangle$
for each $ b \otimes g , a \otimes h \in U$.

The element $\langle u,v \rangle - \langle v,u \rangle$
belongs to every maximal ideal for each $u,v \in U$ 
by Theorem 10 and the consideration above.
If $O$ is semiprime 
the intersection of maximal ideals (i.e. the radical) is zero.
Thus, the form is symmetric.
$\Box$

The following question seems to be of interest. 
Assume that the ground field {\bf k} is
formally real, i.e. zero cannot be presented as 
a sum of squares. 
Under which conditions is the form positive definite. 
Over $\Bbb R$
it is a question of having a Riemannian structure.

\begin{prop} Let {\bf k} be formally real. 
If the form is positive definite  then $H$ is cosemisimple.
\end{prop}
{\bf Proof.}
$0 < \langle 1,1 \rangle =  \Lambda (1) 1$
Thus, $\Lambda (1) \neq 0$ and $H^{\ast}$ is cosemisimple
by the Sweedler criterion. $\Box$

We may also notice that a formally real 
field has zero characteristic.
Thus, $H$ must be as well semisimple  
by the Larson-Radford theorem.

The form allows us to introduce two more features. The map 
$C: U \otimes_O U \otimes_O U \rightarrow {\bf k}$
given by $C(u \otimes v \otimes w)= \langle uv, w \rangle$ 
encodes the multiplication.
There is also  the Nakayama automorphism \cite{nak} defined
by $\langle u,v \rangle = \langle v, \mbox{Nak}(u) \rangle$.
It measures the failure of the form to be symmetric. 
In particular,
Nak is inner if and only if $U$ admits a symmetric associative
bilinear form. 

\section{Connections.}

The study of connections on both principal and 
vector bundles is an important part of modern geometry.
Connections on quantum principal bundles were 
actively studied \cite{haj}.
A Hopf-Galois extension carries also a structure 
of a vector bundle which connections
may have some significance. Moreover, the notion 
of connection on a vector bundle
is less technical than that on a principal bundle 
yet to mention quantum principal bundles.
Throughout this section we prove some interesting 
propositions, involving connections and
Hopf theoretical features of the extensions, 
towards a conjectural  existence and uniqueness 
theorem for a some
kind of connections similar to the existence
and uniqueness of the torsion-free Riemannian 
connection on the tangent
bundle of a Riemannian manifold.

We study  bilinear pairings 
$ \nabla : {\mathcal L} \times U \longrightarrow U$ 
where ${\mathcal L}$ is the Lie algebra 
of derivations of $O$ (vector fields
on the spectrum).
We denote the result of the pairing by 
$\nabla_X u$ with $X$ being a vector field.
 $\nabla$ is called {\em a connection} if 
$$\nabla_X  (au)  =  X(a)u+ a \nabla_Xu \mbox{ and }
\nabla_{aX}u=a \nabla_X  u.$$
for each $a \in O$, $X \in {\mathcal L}(O)$, and $u \in U$.
A connection $\nabla$ is called {\em a Frobenius connection} if
$$X \langle u, v \rangle = 
\langle \nabla_X  u, v \rangle + \langle u , 
\nabla_X  v \rangle.$$
A connection $\nabla$ is called {\em multiplicative} if
$$\nabla_X  (u v)  =  ( \nabla_X  u ) v  +  u  \nabla_X  v.$$
A connection $\nabla$ is called {\em equivariant} if
$$\nabla_X  u_0 \otimes u_1 = 
( \nabla_X  u)_0 \otimes ( \nabla_X  u)_1$$
A connection $\nabla$ is called {\em Nakayama} if
$$\mbox{Nak}^{-1}( \nabla_X  \mbox{Nak}(u))=  \nabla_X  u.$$

Roughly speaking Frobenius connection 
is the one along which the form is covariant constant.
Indeed, being Frobenius is equivalent to 
$\nabla_X \langle , \rangle =0$.
A multiplicative connection may be thought of as a way 
to extent derivation of $O$
to derivations of $U$.  The properties we have 
just defined are not independent.

\begin{prop}
A multiplicative equivariant connection is Frobenius. 
A Frobenius connection 
is Nakayama.
\end{prop}
{\bf Proof.} For each $u,v \in U$ we have 
$X ( \langle u, v \rangle ) = 
\nabla_X (u_0 v_0 \Lambda (u_1v_1))=$
$$\nabla_X (u_0) v_0 \Lambda (u_1v_1) + 
u_0 \nabla_X (v_0) \Lambda (u_1v_1) +
u_0 v_0 X( \Lambda (u_1v_1) )=$$
$$=( \nabla_X u)_0 v_0 \Lambda 
( ( \nabla_X u)_1v_1) + u_0 ( \nabla_X v_0) \Lambda (u_1
( \nabla_Xv)_1) = \langle  \nabla_X u, v 
\rangle +  \langle u,  \nabla_X v \rangle $$
which proves the first claim. 
To show the second one it suffices to notice the following
sequence of equalities for all 
$u,v \in U$ and $X \in {\mathcal L}$
$$\langle \nabla_X u,v \rangle 
=X(\langle u,v \rangle ) - \langle u, \nabla_Xv \rangle=
X(\langle v, \mbox{Nak}(u) \rangle ) - 
\langle  \nabla_Xv, \mbox{Nak}(u) \rangle=$$
$$=\langle  v, \nabla_X \mbox{Nak}(u) \rangle
=\langle  \mbox{Nak}^{-1}(\nabla_X 
\mbox{Nak}(u)),v \rangle \ \ \ \ \ \ \ \ \ \ \ \Box$$

The next three propositions deal with 
the issue of existence of a connection
with certain properties.

\begin{prop}
A connection always exists.
\end{prop}
{\bf Proof.}
$U$ is direct summand of a free $O$-module
$O^n$. We identify $U$ with the submodule of $O^n$. 
Let $\pi : O^n \rightarrow U$ be
a projection. Let $e_i$ be a base of $O^n$. 
The free module admits a trivial connection
making $e_i$ covariant constants. We want 
to restrict this connection to $U$.
Given $u=u^ie_i \in U$, we define $\nabla_Xu = \pi (X(u^i)e_i)$.
The map $\nabla$ is obviously a connection: for instance, 
$\nabla_X au=\pi (X(au^i)e_i)= \pi (X(a)u^i e_i)+ 
\pi (aX(u^i)e_i)= X(a)\pi (u^ie_i)+ a \pi (X(u^i)e_i)
=X(a)u+ a\nabla_Xu$ for $X \in {\mathcal L}, 
a \in O, u \in U$. $\ \ \ \Box$

\begin{prop}
If the characteristic of {\bf k} is not 2 and 
there exists a Nakayama connection
then there exists a Frobenius connection.
\end{prop}
{\bf Proof.}
Let  $\nabla$ be a Nakayama connection. 
Let us define $T$ by $\langle T_Xu,v \rangle = X( \langle 
u,v \rangle )- \langle \nabla_Xu,v \rangle - 
\langle u, \nabla_X v \rangle $.
Then $\langle u,T_Xv \rangle = 
\langle T_Xv, \mbox{Nak}(u) \rangle
= X( \langle 
v, \mbox{Nak}(u) \rangle )- 
\langle \nabla_Xv, \mbox{Nak}(u) \rangle
- \langle v, \nabla_X \mbox{Nak} (u) \rangle
=X( \langle u,v \rangle )- \langle u, \nabla_Xv \rangle
- \langle \mbox{Nak}^{-1}( \nabla_X \mbox{Nak}(u)),v \rangle $.

It is easy to show that $\tilde{\nabla}
= \nabla + \frac{1}{2}T$ is a connection:
$\langle \tilde{\nabla}_{aX}bu , v \rangle =
\langle \nabla_{aX}bu , v \rangle +
\frac{1}{2}[aX( \langle bu,v \rangle )
- \langle \nabla_{aX}bu,v \rangle 
- \langle bu, \nabla_{aX} v \rangle]=
ab \langle \nabla_{X}u , v \rangle 
+ aX(b) \langle u, v \rangle +
\frac{1}{2}[aX(b) \langle u, v \rangle+
ab X( \langle u,v \rangle )- aX(b) \langle u, v \rangle
-ab \langle \nabla_{X}u,v \rangle - 
ab \langle u, \nabla_{X} v \rangle]=
 ab \langle \tilde{\nabla}_{X}u , 
v \rangle  + aX(b) \langle u , v \rangle +
\frac{1}{2}ab  \langle T_Xu,v \rangle ).$

$\tilde{\nabla}$ is Frobenius:
$ \langle \tilde{\nabla}_X u, v \rangle 
+ \langle u, \tilde{\nabla}_X v \rangle=
\langle \nabla_X u, v \rangle + 
\frac{1}{2}[X( \langle u,v \rangle )- 
\langle \nabla_Xu,v \rangle - \langle u, \nabla_X v \rangle]+
\langle u, \nabla_X v \rangle +\frac{1}{2}[
X( \langle u,v \rangle )- \langle \nabla_Xu,v \rangle -
 \langle u, \nabla_X v \rangle]. 
\ \ \ \ \ \ \ \Box$

The next lemma describes another trick 
which may be performed with connections.
It may be restated that a certain affine 
hyperplane of $H^\ast$ is acting
on the space of connections.

\begin{lemma} Let $h$ be an element of $H^{\ast}$ such that
$h(1)=1$ and $\nabla$ be a connection. 
The pairing $h \cdot \nabla$
defined by  $(h \cdot \nabla)_Xu= 
h_1 \cdot (\nabla_X Sh_2 \cdot u)$
is a connection.
\end{lemma}
{\bf Proof.} 
It is clear that  $(h \cdot \nabla)_{fX}u = 
h_1 \cdot (\nabla_{fX} Sh_2 \cdot u)
= h_1 \cdot (f \nabla_X Sh_2 \cdot u)
= (h_1 \cdot f)h_2 \cdot  (\nabla_X Sh_3 \cdot u)
=f(h \cdot \nabla)_Xu$ for each 
$X \in {\mathcal L}, u \in U$.

It suffices to check
$ (h \cdot \nabla)_Xfu =h_1 \cdot [\nabla_X (Sh_2 \cdot fu)]=
h_1 \cdot [\nabla_X ((Sh_3 \cdot f)( Sh_2 \cdot u))]=
 h_1 \cdot [\nabla_X  f( Sh_2 \cdot u)] = 
h_1 \cdot [X ( f)( Sh_2 \cdot u)
+ f \nabla_X ( Sh_2 \cdot u)] = 
h_1 \cdot [X ( f)( Sh_2 \cdot u)]
+ h_1 \cdot [f \nabla_X ( Sh_2 \cdot u]=
h(1)X(f)u+f(h \cdot \nabla)_X u = 
X(f)u+f(h \cdot \nabla)_X u. \ \ \ \ \ \Box$

We remind that 
$H$ being cosemisimple 
is equivalent to  $\Lambda (1) \neq 0$
for some integral $\Lambda$. 
Without loss of generality we may assume
that $\Lambda (1)=1$. 
This allows us to integrate connections: 
given $\nabla$, we obtain $\Lambda \cdot \nabla$.
We also remind that in zero characteristic 
a cosemisimple Hopf algebra is involutive
and unimodular (and even semisimple) \cite{lr1,lr2}.

\begin{prop} If $H$ is  cosemisimple
then $\Lambda \cdot \nabla$ is an  equivariant 
connection  for any connection $\nabla$.
If $\nabla$ is multiplicative 
then so is $\Lambda \cdot \nabla$
provided $H$ is unimodular and involutive.
\end{prop}
{\bf Proof.} We utilize \cite[Formula 3]{rad} 
with $a=S^{-1}(h) \in H^\ast$ 
to prove the first statement:
$$h \Lambda_1 \otimes \Lambda_2 
= \Lambda_1 \otimes S^{-1}(h) \Lambda_2.$$
Let us apply $\mbox{Id} \otimes S$ to this:
$$h \Lambda_1 \otimes S(\Lambda_2) = 
\Lambda_1 \otimes S( \Lambda_2)h.$$
Now, $h \cdot [(\Lambda \cdot \nabla)_Xu] = 
h \Lambda_1 \cdot [ \nabla_X (S(\Lambda_2) \cdot u)] =
\Lambda_1 \cdot [ \nabla_X (S(\Lambda_2)h \cdot u)] = 
(\Lambda \cdot \nabla)_X h \cdot u$
which is reformulation of the equivariance  
condition in the language of $H^\ast$-action.

$H$ being unimodular and involutive is equivalent 
to $\Lambda$ being cocommutative.
It  easily implies  that 
$\Lambda_1 \otimes \Lambda_2 \otimes 
\Lambda_3 \otimes \Lambda_4 =
\Lambda_4 \otimes \Lambda_1 \otimes 
\Lambda_2 \otimes \Lambda_3.$
If $f \in O$ then  $(\Lambda \cdot \nabla)_X(f) = 
\Lambda_1 \cdot 
\nabla_X(S \Lambda_2 \cdot f)
= \Lambda (1) X(f)=X(f)$. 
$H^\ast$ is involutive since so is $H$. 
We are ready to check that $\Lambda \cdot \nabla$
is multiplicative provided so is $\nabla$:
$$(\Lambda \cdot \nabla)_X(uv)= \Lambda_1 
\cdot \nabla_X(S \Lambda_2 \cdot uv)=
\Lambda_1 \cdot \nabla_X[(S \Lambda_3 \cdot u) 
(S \Lambda_2 \cdot v)]=$$
$$=[\Lambda_1 \cdot \nabla_X(S \Lambda_4 \cdot u)] 
(\Lambda_2 S \Lambda_3 \cdot v)
+\Lambda_1 S \Lambda_4 \cdot u [ \Lambda_2 \cdot 
\nabla_X (S \Lambda_3 \cdot v)]=$$
$$=\Lambda_1 \cdot [\nabla_X(S \Lambda_2 \cdot u)]  v
+\Lambda_2 S^{-1} \Lambda_1 \cdot u [ \Lambda_3 
\cdot \nabla_X (S \Lambda_4 \cdot v)]=$$
$$= [(\Lambda \cdot \nabla)_X u]  v
+u (\Lambda \cdot \nabla)_X v. 
\ \ \ \ \ \  \ \ \ \ \Box$$

Now we would like to address an issue 
of uniqueness of  connection.
Let us assume that we are given  
two connections $\nabla^1$
and $\nabla^2$. Let $D= \nabla^1- \nabla^2$. 
It is straightforward to see
that $D_X$ is an endomorphism of the $O$-module $U$.
Thus, two connections differ by an element of 
$\Hom_O ({\mathcal L}, \mbox{End}_OU)$.  According to
\cite[8.3.3]{mon}, $ \mbox{End}_OU \cong U \# H^\ast$, 
the smash product, as algebras.
In particular, $\Hom_O ({\mathcal L}, \mbox{End}_OU)
\cong \Hom_O ({\mathcal L}, U \# H^\ast )
\cong \Hom_O ({\mathcal L}, U)^n$ where $n$ 
is the dimension of $H$.
We have just proved the following proposition.

\begin{prop} The space of connections is a 
$\Hom_O ({\mathcal L}, U)^n$-torsor
(i.e. a space with free transitive action of  
$\Hom_O ({\mathcal L}, U)^n$).
 \end{prop}

If the connections  $\nabla^1$
and $\nabla^2$ are equivariant then 
the difference  $D= \nabla^1- \nabla^2$
is an element of  $\Hom_O ({\mathcal L}, 
\mbox{End}_{O \otimes H^\ast}U)$
and the space of equivariant connections 
is a torsor over this group.
It allows a nice interpretation 
if $U$ is cleft: let us assume
$U \cong O_\sigma [H]$ is given.
Now, $\mbox{End}_{O \otimes H^\ast}U 
\cong \mbox{Comod} (H,U)$.
Right $H$-comodule maps from $H$ to 
$U$ (or integrals in the
terminology of \cite{doi,dot}) 
were actively studied. We think
that the fact we have just noticed 
is worth writing as a proposition.
\begin{prop} 
Given an isomorphism $U \cong O_\sigma [H]$, the space
of equivariant connections becomes a
 $\Hom_O ({\mathcal L}, \mbox{Comod}(H, U))$-torsor.
\end{prop}

Adding other special  properties 
in this spirit will further refine
the space of connections.
For example, if  $\nabla^1$
and $\nabla^2$ are multiplicative connections 
then $\nabla^1_X- \nabla^2_X$ is 
a differentiation of $U$ for each $X$.
Such process may eventually lead 
to an interesting uniqueness theorem.

The last proposition is just 
a reformulation of Frobenius property involving
the map $C$ defined in Section 5.

\begin{prop} We assume that $\nabla$ is Frobenius.
Given $X$, $\nabla_X$ is a derivation of $U$ if and only
if $\nabla_X \cdot C=0$
\end{prop}
{\bf Proof.}
Being a derivative means the equality
$$\langle \nabla_X (uv),w \rangle = 
\langle \nabla_X (u)v,w \rangle +
\langle u \nabla_X (v),w \rangle$$
for each $u,v,w \in U$. Since 
the connection is Frobenius it may be rewritten
$$X( C (u \otimes v  \otimes w)) - 
\langle uv, \nabla_X w \rangle = 
\langle \nabla_X (u)v,w \rangle +
\langle u \nabla_X (v),w \rangle$$
By the definition,
$$\nabla_X \cdot C(u  \otimes v  \otimes w) 
= X(C(u  \otimes v  \otimes w))
-C( \nabla_Xu  \otimes v  \otimes w) - 
C(u  \otimes \nabla_Xv  \otimes w)-
C(u  \otimes v  \otimes  \nabla_Xw)$$
This proves the proposition. $\Box$

\begin{cor} A Frobenius connection 
is multiplicative if and only if
$\nabla \cdot C=0$.
\end{cor}

It may be also interesting 
to investigate a significance of  flatness.
A connection is flat if the curvature form 
$R_{X,Y}u= [ \nabla_X , \nabla_Y ]u - 
\nabla_{[X,Y]}u$ is zero
which is equivalent to $\nabla$ defining 
a representation of Lie algebra ${\mathcal L}$ on $U$.

\section{Miyashita-Ulbrich action.}
Miyashita-Ulbrich action was introduced in \cite{dot}.
If $U \subseteq O$ is 
an $H$-Galois extension with not necessarily central $O$
we can define the Miyashita-Ulbrich action 
for any algebra morphism $\alpha : U \longrightarrow B$. 
We denote $B^O$ and $B^U$ the centralizer
of $\alpha (O)$ and $\alpha (U)$ in $B$. 
Miyashita-Ulbrich action is 
a right $H$-action on $B^O$ defined
by the formula 
$c \stackrel{\alpha}{\leftarrow} x = 
\sum \alpha (a_i) c \alpha (b_i)$
for any $x \in H$ so that $\sum a_i \otimes b_i
= \mbox{can}^{-1} (1 \otimes x)$.
It has a property that for $c \in B^O$ and $b \in U$
the formula $c \alpha (b) = \alpha (b_0) (c  \stackrel{\alpha}
{\leftarrow} b_1)$ holds. Furthermore, 
the invariants of the Miyashita-Ulbrich
action is precisely $B^U$ \cite{dot}.

Keeping in mind that $U \supseteq O$ 
is an extension with central
invariants we consider  the identity 
map from $U$ to $U$ as $\alpha$.
This defines a right $H$-action on $U$ 
such that the center of
$U_{\chi}$ is the subalgebra of invariants. 
We use a shorthand notation
$\leftarrow$ rather than 
$\stackrel{\small Id}{\longleftarrow}$ for this action.
Using antipode one
can get a left Miyashita-Ulbrich as well. 
By doing so one obtains a left action
of the quantum double $D(H)$ 
in the case of cocommutative $H$ \cite{coh}.

Miyashita-Ulbrich action seems 
to be an important piece of structure in our set-up. 
We  think
of $U_{\chi}$ as a deformation of  right $H$-module. 
One may also think of $U_{\chi}$ 
as a deformation
of an algebra structure: it was done in  \cite{rum}. 

 Let $E_{\chi}= \mbox{End}_{\bf k}U_{\chi}$. 
Then $E_{\chi}$
carries a structure of $H$-$H$-bimodule by 
$h \phi h^{\prime} (u) = \phi (u \leftarrow h) 
\leftarrow h^{\prime}$. We are interested in  
the first Hochschild cohomology HH$^1(H,E_\chi )$
because normalized first Hochschild cocycles 
with coefficients in $E_\chi$
 constitute infinitesimals of deformations 
of module structure.
Indeed, let $D={\bf k}[ \epsilon] / (\epsilon^2)$ 
be the ring of dual numbers. An infinitesimal
deformation is a right $H$-module structure 
$U_\chi [\epsilon] \otimes_DH \rightarrow U_\chi [\epsilon]$
of the form:
$$a \otimes h \mapsto \ \ \ 
a \leftarrow h + \mu (a,h) \epsilon$$
for $a \in U_\chi, h \in H$. 
The following is the associativity condition:
$$a \leftarrow hh^\prime + \mu (a,hh^\prime ) 
\epsilon =a \leftarrow h \leftarrow h^\prime + 
\mu (a,h) \leftarrow h^\prime \epsilon +  
\mu (a \leftarrow h, h^\prime ) \epsilon .$$
It may be rewritten  as 
$\mu (a,hh^\prime ) =\mu (a,h) \leftarrow h^\prime + 
 \mu (a \leftarrow h, h^\prime )$
which is the Hochschild 1-cocycle condition for the map 
$\tilde{\mu}: H \rightarrow E_\chi$
given by $\tilde{\mu} (h) (a)= \mu (a,h)$. 
The unitary condition $\mu (a,1)=0$
corresponds to the normalization condition on the cocycle.

We should point out that this consideration 
never uses specific features of the situation:
it works for any algebra $H$ and a family of 
right $H$-modules $U_\chi$.
The following theorem is an adaptation 
of \cite[theorem 3.2]{ger} for module deformations.
One can probably prove a stronger version of 
this theorem utilizing other ideas
of \cite{ger}.

\begin{theorem}
We assume {\em HH}$^1(H,E_{\chi})=0$ for some $\chi$. 
If $C$ is a curve in $\mbox{\em Spec} \, O$
smooth at the point $\chi$ 
then there exists
an open neighborhood 
$W$ of $\chi$ in $C$ such
that 
$U_{\eta} \otimes L \cong U_{\chi} \otimes L$ 
as $H$-modules under Miyashita-Ulbrich action 
for each $\eta \in W$
and some field extension $L \supseteq {\bf k}$.
\end{theorem}
{\bf Proof.}
Let $t$ be a local parameter on $C$ at the point $\chi$.
It gives us a generic point $\beta (t)$ 
of $C$ such that $\beta (0)= \chi$.
The generic point may be thought of 
as an element of  $\hom (O, K[[t]])$
where $K$ is the field $O/\chi$.

By Lemma 5, $U \otimes_O K[[t]]$ 
is an $H$-Galois extension with
central invariants. Therefore, 
it experiences the Miyashita-Ulbrich action
which may be written as 
a formal deformation of that of $U_\chi$:
$$\phi : U_{\chi}[[t]] \otimes_{{\bf k}[[t]]} 
H \longrightarrow U_{\chi}[[t]].$$
$\phi$ may be given through a bunch of 
$\phi_n : U_{\chi} \times H \rightarrow U_{\chi}$
such that $\phi (a \otimes  h) = \sum_n \phi_n(a ,h)t^n$.
Clearly, $\phi_0 (a \otimes h)= a \leftarrow h$.  
The non-zero $\phi_n$ with the smallest
positive $n$ is a first normalized Hochschild 
1-cocycle with coefficients in $E_\chi$
by the argument similar to the one about infinitesimals.

But every cocycle is a coboundary! Thus, 
there exists $\Theta \in E_\chi$ such that
$\phi_n (a,h) = 
\Theta (a \leftarrow h)- \Theta (a) \leftarrow h$. 
We consider an $H$-module isomorphism 
$F_n: U_\chi[[t]] \rightarrow U_\chi [[t]]$ given by
$F_n (a)=a+ \Theta (a)t^n$. It is clear that 
$F_n^{-1}(a)=a- \Theta (a)t^n + {\frak o}(t^n)$.
Let us compute
$$F_n^{-1} (F_n(a)h)= F_n^{-1} (a \leftarrow h +
\Theta (a)\leftarrow ht^n + \phi_n (a,h)t^n +
{\frak o}(t^n))= $$
$$a \leftarrow h + [-\Theta (a \leftarrow h) + 
\Theta (a)\leftarrow h + \phi_n (a,h)]t^n +
{\frak o}(t^n))= a \leftarrow h + {\frak o}(t^n).$$
Thus, $F_n$ ``kills''  the $t^n$-term of the deformation. 
We may go on defining $F_{n+1},
F_{n+2}, \ldots$ in the similar fashion. The product 
$F= \cdots F_{n+1} \circ F_n$, though
infinite, is well-defined because the coefficient 
at $t^k$ is computed in the finite
product $F_k \circ \cdots \circ F_n$. It is clear that
$$F^{-1}(F(a)h) = a \leftarrow h$$
and, therefore, $F$ performs an $H$-module 
isomorphism between $U_\chi [[t]]$ and
the trivial deformation $U_\chi \otimes_K K[[t]]$. 
But $U_\chi [[t]]$ may be specialized
to $U_\eta$ for $\eta$ from some open subset $W$ of $C$. 
Thus, $U_\eta \otimes_{O/\eta} K((t))$ must be isomorphic
to $U_\chi \otimes_K K((t))$ as right $H$-modules.
$\Box$

Theorem 26 gives a fine tool for calculations. 
Let $\frak g$ be a semisimple classical Lie
algebra of dimension $n$ and rank $r$.
We assume that $\bf k$ is algebraically 
closed of good characteristic 
(i.e. the quotient of the root lattice 
by any sublattice generated by 
a root subsystem has no $p$-torsion). 
The list of bad primes is
2 for $B_n$, $C_n$, $D_n$, 2 and 3 for 
$E_6$, $E_7$, $F_4$, $G_2$,
and 2, 3, and 5 for $E_8$. 
One should consult \cite[3.2]{rum} or 
\cite{fri} for the definition of reduced
enveloping algebras. Roughly speaking, 
they are algebras $U_\chi$ as $U=U({\frak g})$
and $O=O({\frak g}) \cong 
{\mathcal O}({\frak g}^{\ast (1)})$. We consider a line
$\alpha \chi$ with $\alpha \in {\bf k}$ 
and regular semisimple $\chi$.
$U_0$ is just restricted enveloping algebra. 
By \cite[Corollary 3.6]{fri} 
$U_{\alpha \chi}$ for $\alpha \neq 0$ 
is a direct sum of $p^r$ copies
of the algebra of $p^{n-r} \times p^{n-r}$-matrices. 
Since the dimension of the center
does not change with a field extension 
(by an elementary argument
as in \cite[Theorem 26]{rum}) 
we obtain the following corollary 
of Theorem 26.

\begin{cor} Let $\frak g$ be a classical semisimple Lie
algebra. We assume $p$ is good.
If 
{\em HH}$^1(u({\frak g}), \mbox{End}_{\bf k}u({\frak g}))=0$ 
where $u({\frak g})$ is the module over itself
under the right adjoint action then the center of 
$u({\frak g})$ has dimension $p^r$.
\end{cor}

It follows from the results of \cite{kac} 
that the dimension of the center is at least $p^r$. 
To be precise one needs
slightly stronger restriction on $p$: 
$p$ has to be good and should not divide 
$n+1$ if ${\frak g}$
has $A_n$ as a direct summand. 
To the best of our 
knowledge, it is unknown whether 
the dimension is precisely $p^r$ even for $sl_2$.
We don't know how to compute 
HH$^1(u({\frak g}), \mbox{End}_{\bf k}u({\frak g}))$
but it looks like a possible way to put hands on the center. 
Not only could it settle the question
about the center of $u({\frak g})$ but 
it could also help to understand  
the centers of $U_\chi$ for 
nilpotent $\chi$. Jacobson-Morozov theorem provides 
that $\chi$ and $\alpha \chi$ are conjugate under
coadjoint action for $\alpha \neq 0$ which implies 
$U_{\alpha \chi} \cong U_\chi$.
We can use Theorem 26 to the curve ${\bf k} \chi$
at 0.

\begin{cor} Let $\frak g$ be a classical semisimple
Lie algebra with the Coxeter number $h$.
If $p > 3h-1$ and
{\em HH}$^1(u({\frak g}), \mbox{End}_{\bf k}u({\frak g}))=0$ 
where $u({\frak g})$ is the module over itself
under the right adjoint action then 
the center of $U_\chi$ has dimension $p^r$ 
for every nilpotent $\chi$.
\end{cor}

\section{$H$-torsors on schema.}

{\em An $H$-extension} of a {\bf k}-ringed space 
$(X, {\mathcal R})$ (i.e.
a topological ring with a sheaf of algebras) 
is a sheaf of $H$-comodule
algebras $\mathcal U$
on $X$ (i.e. $H$ coacts on ${\mathcal U}(V)$ 
for each open set $V$ and restriction
maps are $H$-equivariant) such that $\mathcal R$ 
is a subsheaf of $\mathcal U$ and
${\mathcal R}(V) ={\mathcal U}(V)^H$ 
for each open set $V$. 

This definition is not very helpful. 
The problem is that we cannot glue
such objects from local data. 
To check that we have an $H$-extension
we must look at each open set.  For instance, 
if $X$ is an affine scheme
and $\mathcal R$ is a sheaf of algebraic 
functions $\mathcal O$ then 
a Hopf-Galois extension of the algebra 
${\mathcal O}(X)$ does not
necessarily give a Hopf-Galois extension of 
$(X,{\mathcal O})$.
However, if we restrict our attention
to $H$-Galois extensions with central invariants 
then it suffices
to work with open sets from some affine covering.

We are interested only in a scheme $\mathcal S$ over 
{\bf k} and the sheaf of regular functions
$\mathcal O$.
By an {\em $H$-torsor on the scheme $\mathcal S$} 
we understand an $H$-extension
$\mathcal U$ of $({\mathcal S,O})$ with central invariants 
which is locally $H$-Galois.
The latter means that there exists an affine 
covering $\coprod_i V_i$ of $\mathcal S$
in the flat topology such that ${\mathcal U}(V_i)$ 
is a Hopf-Galois extension of ${\mathcal O}(V_i)$
for each $i$. We should point out that 
the extension 
${\mathcal U}(V_i) \supseteq {\mathcal O}(V_i)$
has central invariants automatically.

One can use the usual gluing procedure to paste 
$H$-torsors. $H$ torsors on an affine scheme
is the same as $H$-Galois extensions of the algebra 
of global algebraic functions. {\em An $H$-structure}
on a scheme $\mathcal S$ is the following data:
an affine covering  $\coprod_i V_i$, a collection 
of $H$-Galois extensions $U_i \supseteq {\mathcal O}(V_i)$, 
and a collection of isomorphisms
$\phi_{ij}:{\mathcal U}(V_i) \otimes_{\mo (V_i)} 
{\mathcal O}(V_i \times_{\mathcal S}V_j)
\rightarrow {\mathcal U}(V_j) 
\otimes_{\mo (V_j)} {\mathcal O}(V_i \times_{\mathcal S}V_j)$
which satisfy the cocycle condition on triple intersections 
and give the identity map
upon restriction to $\mo (V_i \times_\sss V_j)$.

Similarly to Theorem 11 an $H$-structure determines 
the $H$-torsor.
However, an $H$-torsor is not necessarily a sheaf
of $H$-Galois extensions. The problem is that
the extension of sections on a bad open
subset can fail to be Galois.
For each scheme $\sss$ there is a natural map 
$\psi_{\sss} : \sss \rightarrow \mbox{Spec}\, \mo (\sss )$.
We call the scheme {\em good} if  the natural map 
$\psi_{\sss}$ is faithfully flat.
Both affine and projective schemes are good but 
there are  schemes which are not.
If, for instance, $\sss$ is the complement of 
a point in ${\Bbb A}^2$ then 
$\mbox{Spec}\, \mo (\sss ) \cong {\Bbb A}^2$ 
and $\psi_{\sss}$ is the natural embedding which is not
surjective and, therefore, not faithfully flat.

\begin{theorem}
If $V$ is a good open subscheme of $\sss$ and 
$\uu$ is an $H$-torsor on $\sss$  then
the extension $\uu (V) \supseteq \mo (V)$ is $H$-Galois.
\end{theorem}
{\bf Proof.}
Let $\coprod_i V_i$ is an affine open covering
of $\sss$ such that $\uu (V_i) \supseteq \mo (V_i)$
is $H$-Galois for each $i$.
Let $\coprod_j W_{ij}$ be an affine open covering
of $V \times_\sss V_i$. 
The extensions of sections on each $W_{ij}$ is $H$-Galois.
Now the composition
$\coprod_{ij} W_{ij} \rightarrow V \rightarrow 
\mbox{Spec}\, \mo (V)$ is faithfully flat 
since so are both maps. Theorem 11, applied to
the composition, proves Theorem 26. $\ \Box$

We need the notion of pull-back of $H$-torsors
to continue the discussion. 
Given a map of schemes 
$f:{\mathcal S}^\prime \rightarrow {\mathcal S}$ 
and an $H$-torsor $\mathcal U$
on $\mathcal S$, there exists 
{\em a pull-back $H$-torsor} $f^\ast {\mathcal U}$. 
It is constructed by  choosing 
compatible coverings as in Theorem 29 and applying
tensor products (Lemma 5).
We obtain an $H$-structure which determines
the $H$-torsor. 
The definition is independent of the choice of covering:
given two coverings, one can use their common refinement
to construct the pull-back. This construction also
gives a canonical isomorphism between the two $H$-torsors
constructed by the two coverings.

A natural question is to classify $H$-torsors 
on a scheme $\mathcal S$. As one could expect this problem
leads to theory of stacks \cite{del,vis}. 
Let us build a category $\ca$. 
The objects of $\ca$ are pairs $(X,\uu )$
where $X$ is a scheme over {\bf k} and $\uu$ 
is an $H$-torsor on $X \times \sss$.
Morphisms $F:(X^\prime ,\uu^\prime)
\rightarrow (X, \uu )$ 
are maps $f:X^\prime \rightarrow X$ such that
$(\ii_{\uu} \times f)^\ast (\uu ) \cong \uu^\prime$. 

\begin{theorem} The category $\ca$ together with 
the forgetful functor to $\pi$ the category of scheme
(i.e. $\pi (X, \uu)= X$) is a stack in the flat topology
over the category of schemes over {\bf k}.
\end{theorem}
{\bf Proof.} We go over the four defining 
properties of stacks as defined in \cite{del,vis}. 
The first two properties determine a groupoid
over the category of schemes, the others provide
that a groupoid is a stack.

1. If $f:X \rightarrow Y$ is a map of schemes and 
$(Y, \uu)$ is an object of $\ca$ then $f$
is lifted to a map in $\ca$ from 
$(X, (\ii_{\sss} \times f)^\ast (\uu ))$ to $(Y, \uu )$.

2. If $c:X \rightarrow Y$ and $b:Y \rightarrow Z$ are maps 
of schemes and $a=b \circ c$ and 
and $(X, \uu )$, $(Y, \uu_1)$, and $(Z, \uu_2)$ are objects 
in $\ca$ such that $a$ and $b$
define morphisms in $\ca$ between the corresponding objects 
then $c$ also defines 
the unique (uniqueness follows from the definition 
of the morphism in $\ca$) morphism
from $(X, \uu )$ to $(Y, \uu_1 )$ because 
$(\ii_{\sss} \times c)^\ast
 \circ (\ii_{\sss} \times b)^\ast (\uu_3) 
\cong (\ii_{\sss} \times b \circ c)^\ast (\uu_3)= 
(\ii_{\sss} \times a)^\ast (\uu_3)$. 
This holds because the pull-back is defined by tensor products
and tensor product is associative.

3. For each scheme $X$ and each two objects
$(X, \uu_1)$ and $(X,\uu_2)$ of $\ca$
the functor $I$ from the category of schemes over $X$
to the category of sets such that 
$I( Y \stackrel{f}{\rightarrow}X)$ is a set
of isomorphisms in $\ca$ 
between $f^\ast (\uu_1)$ and $f^\ast (\uu_2)$
must be a sheaf in the flat topology.
Indeed, let $\coprod_i g_i: \coprod_i Y_i \rightarrow Y$
be a covering. We consider a bunch of isomorphisms
$$\Omega_i: (fg_i)^\ast (\uu_1) \rightarrow
(fg_i)^\ast (\uu_2)$$
which agree on double intersections. 
They are given by isomorphisms 
$\omega_i:Y_i \rightarrow Y_i$ such that
$fg_i=fg_i \omega_i$. Since $\omega_i$ agree 
on double intersections they can be glued into
a map $\omega :Y \rightarrow Y$ such that $f=f\omega$.
We have to check that $\omega$ determines a morphism 
in $\ca$ from $f^\ast (\uu_1)$ to $f^\ast (\uu_2)$ 
which is equivalent to the fact that 
$\omega^\ast (f^\ast (\uu_2)) \cong f^\ast (\uu_1)$.
This condition is local but we know that locally
${\omega_i}^\ast ((fg_i)^\ast (\uu_2)) \cong (fg_i)^\ast (\uu_1)$.

4. The last property is that an $H$-structure determines
an $H$-torsor.
$\ \ \ \Box$

\end{document}